\documentclass[12pt, reqno]{amsart}
\usepackage[arrow,matrix,curve]{xy}

\newtheorem{theorem}{Theorem}[section]
\newtheorem{proposition}[theorem]{Proposition}
\newtheorem{lemma}[theorem]{Lemma}
\newtheorem{corollary}[theorem]{Corollary}

\theoremstyle{definition}
\newtheorem{definition}[theorem]{Definition}
\newtheorem{example}[theorem]{Example}
\newtheorem{remark}[theorem]{Remark}

\newcommand{\R}{\mathbb R}
\newcommand{\Z}{\mathbb Z}
\newcommand{\T}{\mathbb T}

\DeclareMathOperator{\sys}{{\rm sys}}
\DeclareMathOperator{\stsys}{{\rm stsys}}
\DeclareMathOperator{\pisys}{\sys\pi} 
\def\phisys {\phi{\rm sys}}
\DeclareMathOperator{\vol}{{\rm vol}} 

\def\area{{\rm area}} \def\length{{\rm length}}

\def\AJ {{\mathcal A}} \def\XX {{\overline{{X}}}}

\def\g {{\it {g}}}  

\def\ie {{\it i.e.\ }} 
\def\eg {{\it e.g.\ }} 
\def\cf {\hbox{\it cf.\ }}

\def\short {{1-Lipschitz\ }} 

\def\Pr {{\rm Pr}}

\numberwithin{equation}{section}
\numberwithin{figure}{section}

\begin{document}

\author[Ivanov]{Sergei V. Ivanov$^\dagger$} \address{Steklov
Math. Institute, Fontanka 27, RU-191011 St. Petersburg, Russia}
\email{svivanov@pdmi.ras.ru} 
\thanks{$^\dagger$Supported by grants
CRDF RM1-2381-ST-02 and RFBR 02-01-00090}

\author[M.~Katz]{Mikhail G. Katz$^{*}$} \address{Department of
Mathematics and Statistics, Bar Ilan University, Ramat Gan 52900
Israel} \email{katzmik@math.biu.ac.il} \thanks{$^{*}$Supported by the
Israel Science Foundation (grants no.\ 620/00-10.0 and 84/03)}

\title[Generalized degree and Loewner-type inequalities] {Generalized
degree and optimal Loewner-type inequalities$^{+}$}
\thanks{$^+$Published in {\em Israel J. Math.} \textbf{141} (2004),
221-233.}


\begin{abstract}
We generalize optimal inequalities of C.~Loewner and M.~Gromov, by
proving lower bounds for the total volume in terms of the homotopy
systole and the stable systole.  Our main tool is the construction of
an area-decreasing map to the Jacobi torus, streamlining and
generalizing the construction of the first author in collaboration
with D. Burago.  It turns out that one can successfully combine this
construction with the coarea formula, to yield new optimal
inequalities.
\end{abstract}


\maketitle

\tableofcontents

\section{Loewner's and Gromov's optimal inequalities}

Over half a century ago, C. Loewner proved that the least length
$\sys_1(\T^2)$ of a noncontractible loop on a Riemannian 2-torus
$\T^2$ satisfies the optimal inequality
\begin{equation}
\label{11b}
\sys_1 \left( \T^2 \right) ^2 \leq \gamma_2\; \area(\T^2),
\end{equation}
where $\gamma_2=\frac{2}{\sqrt{3}}$, \cf \cite{Pu, CK}.

An optimal generalisation of Loewner's inequality is due to M. Gromov
\cite[pp.~259-260]{Gr3} (\cf \cite[inequality (5.14)]{CK}) based on
the techniques of D. Burago and S. Ivanov \cite{BI94, BI95}.  We give
below a slight generalisation of Gromov's statement.

\begin{definition}
Given a map $f:X\to Y$ between closed manifolds of the same dimension,
we denote by $\deg(f)$, either the algebraic degree of~$f$ when both
manifolds are orientable, or the absolute degree, otherwise.
\end{definition}

We denote by $\AJ_X$ the Abel-Jacobi map of $X$, \cf
formula~\eqref{21b}; by~$\gamma_n$, the Hermite constant, \cf
formula~\eqref{49}; and by $\stsys_1(g)$, the stable 1-systole of a
metric $g$, \cf formula~\eqref{41b}.

\begin{theorem}[M. Gromov]
\label{11a}
Let $X^n$ be a compact manifold of equal dimension and first Betti
number: $\dim(X)=b_1(X)=n$.  Then every metric $g$ on $X$ satisfies
the following optimal inequality:
\begin{equation}
\label{10}
\deg(\AJ_X) \stsys_1(g)^n \leq \left( \gamma_n
\right)^{\frac{n}{2}} \vol_n(g).
\end{equation}
The boundary case of equality in inequality~\eqref{10} is attained by
flat tori whose group of deck transformations is a critical lattice in
$\R^n$.
\end{theorem}
Note that the inequality is nonvacuous in the orientable case, only if
the cuplength of $X$ is $n$, \ie the Abel-Jacobi map $\AJ_X$ is of
nonzero algebraic degree.  Recall that a critical lattice (\ie one
attaining the value of the Hermite constant~\eqref{49}) is in
particular extremal, and being extremal for a lattice implies
perfection and eutaxy \cite{Bar}.

Having presented the results of Loewner and Gromov, we now turn to
their generalisations, described in Sections \ref{onetwo} and
\ref{two}.

Our main tool is the construction of an area-decreasing map to the
Jacobi torus, streamlining and generalizing the construction of
\cite{BI94}, \cf \cite{Gr3}.  An important theme of the present work
is the observation that one can successfully combine this construction
with the coarea formula, to yield new optimal inequalities.
Historical remarks and a discussion of the related systolic literature
can be found at the end of the next section.

\section{Generalized degree and first theorem}
\label{onetwo}
\begin{figure}
\[
\xymatrix { & \frac{\deg(\AJ_X) \stsys_1(g)^b} {\vol_n(g)} \leq
(\gamma_b)^{\frac{b}{2}} \\ Loewner:\; \frac{\sys_1(g)^2} { \area(g) }
\leq \gamma_2 \ar[ur]^{Ivanov}_{Katz} \ar[dr]^{Bangert}_{Katz} & \\ &
\frac{\stsys_1(g) \sys_{n-1}(g)}{\vol_n(g)} \leq \gamma_b' }
\]
\caption{\textsf{Two generalisations of Loewner's theorem, \cf
\eqref{ik} and \eqref{22}}}
\label{fig1}
\end{figure}

Let $(X,g)$ be a Riemannian manifold.  Let $n=\dim(X)$ and $b=b_1(X)$.
Loewner's inequality can be generalized by optimal inequalities in two
different ways, as illustrated in Figure~\ref{fig1}, where the
constants $\gamma_n$ and $\gamma_n'$ are defined in
Section~\ref{2two}, while
\begin{equation}
\label{21b}
\AJ_X:X\to \T^b
\end{equation}
is the Abel-Jacobi map (\cite{Li, Gr2}, \cf \cite[(4.3)]{BanK2})
inducing isomorphism in 1-dimensional cohomology.

The map $\AJ_X: X\to T^b$ lifts to a proper map 
\begin{equation}
\label{29}
\overline{\AJ}_X:\XX \to \R^b,
\end{equation}
where $\R^b$ is the universal cover of the Jacobi torus.  The
corresponding fibers of $\AJ_X$ and $\overline{\AJ}_X$ (\ie fibers
which project to the same point of~$T^b$) are homeomorphic.  On the
other hand, the nonvanishing of the homology class $[\overline F_X]$
of the typical fiber of $\overline{\AJ}_X$ in $H_{n-b}(\XX)$ is far
weaker than the nonvanishing of the homology class of the typical
fiber of $\AJ_X$ in $H_{n-b}(X)$ (and thus the resulting theorem is
more interesting).  For instance, for the standard nilmanifold of the
Heisenberg group (\cf Remark~\ref{25}) we have $[\overline
F_X]\not=0$, but the homology class of the fiber of $\AJ_X$ is
trivial.  The latter is true whenever $n=b+1$.

In the theorems below, to obtain a nonvacuous inequality, we replace
the condition of nonvanishing degree in Gromov's theorem, by the
nonvanishing of the homology class $[\overline F_X]$ of the lift of
the typical fiber of $\AJ_X$ to the universal abelian cover $\XX$ of
$X$, \cf \cite[p.~101]{Gr1} and equation~\eqref{29} below.  In
particular, we must have $b\leq n$ to obtain a nonvacuous inequality.
Following M. Gromov \cite [p.~101] {Gr1}, we introduce the following
notion of generalized degree.

\begin{definition}
Denote by $\deg(\AJ_X)$ the infimum of $(n-b)$-volumes of integral
cycles representing the class~$[\overline F_X]\in H_{n-b}(\XX)$, if
$X$ is orientable, and infimum of volumes of $\Z_2$-cycles, otherwise.
\end{definition}

This quantity is denoted `$\deg$' in \cite [p.~101] {Gr1}.  When the
dimension and the first Betti number coincide, this quantity is a
topological invariant.  In general, of course, it is not.  Yet it is
remarkable that the nonvanishing of this quantity suffices to produce
a nonvacuous volume lower bound of Theorem~\ref{13}, generalizing
Gromov's theorem~\ref{11a}.

\begin{theorem}
\label{13}
Let $X$ be a compact manifold.  Let $n=\dim(X)$ and $b=b_1(X)$, and
assume $b\geq 1$.  Then every metric $g$ on $X$ satisfies the
inequality
\begin{equation}
\label{ik}
{\deg(\AJ_X) \stsys_1(g)^b} \leq (\gamma_b) ^{\frac{b}{2}} {\vol_n
(g)}.
\end{equation}
\end{theorem}

Here the stable 1-systole $\stsys_1(g)$ is defined in formula~\eqref
{41b}.  In particular, for $n=b+1$ we obtain the following corollary.

\begin{corollary}
\label{14}
Let $X$ be a compact manifold.  Let $b=b_1(X)$.  Assume that
$\dim(X)=b+1$, and $[\overline F_X] \not=0$.  Then every metric $g$ on
$X$ satisfies the following optimal inequality:
\begin{equation}
\label{11}
\stsys_1(g)^b \pisys_1(g) \leq (\gamma_b)^{\frac{b}{2}} \vol_{b+1}(g),
\end{equation}
where $\pisys_1(g)$ denotes the least length of a shortest
noncontractible loop for the metric $g$.
\end{corollary}

\begin{remark}
\label{25}
To give an example where equality is attained in
inequality~\eqref{11}, it suffices to take a Riemannian fibration by
circles of constant length, over a flat torus whose group of deck
tranformations is critical (choose the circles sufficiently short, so
as to realize the value of the invariant $\pisys_1(g)$).  This can be
realized, for instance, by compact quotients of left invariant metrics
on the 3-dimensional Heisenberg group (here $b=2$), \cf \cite[p.~101]
{Gr1}.
\end{remark}

The results of the present paper are further generalized in
\cite{BCIK1, BCIK2}.  A nonsharp version of inequalities \eqref{11}
and \eqref{12} for an arbitrary pair $b\leq n$ was proved using
different techniques in \cite{KKS}, resulting in an inequality with an
extra multiplicative constant $C(n)$ depending on the dimension but
not on the metric.  Note that a different sharp generalisation of
Loewner's inequality was studied in \cite{BanK1}:
\begin{equation}
\label{22}
\stsys_1(g) \sys_{n-1}(g) \leq \gamma_b' \vol_n(g) ,
\end{equation}
where $\gamma_b'$ is the Berg\'e-Martinet constant, \cf \eqref{22'}.
The work \cite{BanK2} studies the boundary case of equality of a
further generalisation of inequality \eqref{22}.

We would like to point out an interesting difference between
inequalities \eqref{11} and \eqref{22}.  Namely, inequalities for the
1-systole, such as \eqref{11}, tend to be satified even by the
ordinary (\ie unstable) systole, \cf \cite[3C$_1$]{Gr2} and
\cite[(3.5)]{CK}, albeit with a nonsharp constant.  Meanwhile,
inequality \eqref{22} for ordinary systoles is definitely violated by
a suitable sequence of metrics, no matter what the constant, \cf
\cite{BabK}.  Volume lower bounds in terms of systoles and the study
of the associated constants for 4-manifolds appear in \cite{Ka}.

The work \cite{CK} surveys other universal (curvature-free) volume
bounds and formulates a number of open questions, including the one of
the existence of such a lower bound in terms of the least length of a
nontrivial closed geodesic (perhaps contractible) on $X$, for any
manifold $X$, \cf \cite{NR,Sa}.

\section{Pu's inequality and generalisations}
\label{two}
To state the next theorem, we need to recall the inequality of P. Pu.
We record here a slight generalisation of the inequality from
\cite{Pu}, see also \cite{Iv} for an alternative proof and
generalisations.  Namely, every surface $(S,\g)$ which is not a
2-sphere satisfies
\begin{equation}
\label{pu3}
\pisys_1(\g)^2\leq \frac{\pi}{2}\area(\g),
\end{equation}
where the boundary case of equality in \eqref{pu3} is attained
precisely when, on the one hand, the surface $S$ is a real projective
plane, and on the other, the metric $\g$ is of constant Gaussian
curvature.

The generalisation follows from Gromov's inequality \eqref{tqi} below
(by comparing the numerical values of the two constants).  Namely,
every aspherical compact surface $(S,\g)$ admits a metric ball
\[
B=B_p\left(\tfrac{1}{2}\pisys_1(\g)\right) \subset S
\]
of radius $\tfrac{1}{2}\pisys_1(\g)$ which satisfies
\cite[Corollary~5.2.B]{Gr1}
\begin{equation}
\label{tqi}
\pisys_1(\g)^2 \leq \frac{4}{3}\area(B).
\end{equation}

Let $S$ be a nonorientable surface, and let $\phi:\pi_1(S)\to \Z_2$ be
an epimorphism from its fundamental group to $\Z_2$, corresponding to
a map $\hat\phi:S\to \R P^2$ of absolute degree $+1$.  We define the
``1-systole relative to $\phi$'', denoted $\phisys_1(\g)$, of a metric
$\g$ on $S$, by minimizing length over loops $\gamma$ which are not in
the kernel of $\phi$, \ie loops whose image under $\phi$ is not
contractible in the projective plane:
\begin{equation}
\label{cps}
\phisys_1(\g)=\min_{\phi([\gamma])\not= 0\in \Z_2}\length(\gamma).
\end{equation}

Does every nonorientable surface $(S,g)$ and map $\hat\phi:S\to \R
P^2$ of absolute degree one, satisfy the following relative version of
Pu's inequality:
\begin{equation}
\label{23}
\phisys_1(\g)^2\leq \tfrac{\pi}{2}\area(\g)\quad ?
\end{equation}
This inequality is related to Gromov's (non-sharp) inequality
$(*)_{inter}$ from \cite[3.C.1]{Gr2}, see also \cite [Theorem~4.41]
{Gr3}.  This question appeared in \cite [conjecture 2.7]{CK}.  Let
\begin{equation}
\label{16}
\sigma_2=\sup_{(S,g)} \frac{\phisys_1(\g)^2}{\area(\g)},
\end{equation}
where the supremum is over all nonorientable surfaces $S$, metrics $g$
on $S$, as well as maps $\phi$ as above.  Thus we ask whether
$\sigma_2=\frac{\pi}{2}$.  Calculating $\sigma_2$ depends on
calculating the filling area of the Riemannian circle, \cf Proposition
\ref{32} below.

\begin{proposition}
\label{32}
We have the following estimate: $\sigma_2\in [\frac{\pi}{2}, 2]$.
\end{proposition}

\begin{proof}
If we open up a surface $S$ as above along a shortest essential loop
$\gamma$ (in the sense that $\phi([\gamma])\not=0$), we obtain a
$\Z_2$-filling $\Sigma$ (possibly nonorientable) of a circle of length
$2\length(\gamma)$.  It is clear that the boundary circle is imbedded
in $\Sigma$ isometrically as a metric space.  Thus it suffices to
prove that that the filling area of a circle of length $2 \pi$ is at
least $\frac{\pi^2}{2}$.  Choose two points on the boundary circle at
distance $\frac{\pi}{2}$ from each other.  Consider the map $\Sigma
\to \R^2$ whose coordinate functions are the distances to these
points.  The map is area-decreasing, and its image contains a square
of area $\frac{\pi^2}{2}$ (encircled by the image of the boundary),
namely the Pythagorean ``diamond'' inside the square $[0,\pi] \times
[0,\pi]$ in the plane.
\end{proof}

\section{Second theorem: the case $\dim(X)=b_1(X)+2$}

\begin{theorem}
\label{21}
Let $X$ be a compact nonorientable manifold.  Let $b=b_1(X)$.  Assume
$\dim(X)=b+2$ and $\iota [\overline F_X]\not=0$, where $\iota:
H_2(\XX,\Z_2)\to H_2(K,\Z_2)$ is the homomorphism induced by a map to
some aspherical space $K$.  Then every metric $g$ on $X$ satisfies the
following inequality:
\begin{equation}
\label{12}
\stsys_1(g)^b \pisys_1(g)^2 \leq \sigma_2 \gamma_b^{\frac{b}{2}}
\vol_{b+2}(g),
\end{equation}
where $\sigma_2$ is the optimal systolic ratio from \eqref{16}.
\end{theorem}

The proof appears at the end of Section~\ref{seven}.

\begin{remark}
If, as conjectured, we have $\sigma_2=\frac{\pi}{2}$, then the
boundary case of equality in inequality~\eqref{12} is attained by
Riemannian submersions over a flat critical torus, with minimal fibers
isometric to a fixed real projective plane with a metric of constant
Gaussian curvature.
\end{remark}

\begin{example}
For $X=\R P^2 \times \T^2$, we obtain the following inequality:
$\stsys_1(g)^2 \pisys_1(g)^2 \leq \sigma_2 \gamma_2 \vol_4(g),$ which
can be thought of as a ``Pu-times-Loewner'' inequality, \cf
\eqref{11b} and \eqref{pu3} (particularly if we prove that
$\sigma_2=\frac{\pi}{2}$).
\end{example}

\section{Lattices, 
Hermite and Berg\'e-Martinet constants}
\label{2two}
Given a lattice $L\subset (B, \|\cdot\|)$ in a Banach space $B$ with
norm~$\|\cdot\|$, denote by $\lambda_1(L)= \lambda_1(L,\|\cdot\|)>0$
the least norm of a nonzero vector in $L$.  Then the Hermite constant
$\gamma_n>0$ is defined by the supremum
\begin{equation}
\label{49}
\sup_{L\subset \R^n}\frac{\lambda_1(L)^n}{\vol(\R^n/L)} =
(\gamma_n)^{\frac{n}{2}} ,
\end{equation}
where the supremum is over all lattices with respect to a Euclidean
norm.  The particular choice of the exponent ${\frac{n}{2}}$ may be
motivated by the {\em linear\/} asymptotic behavior of $\gamma_n$ as a
function of $n\to\infty$, \cf \cite[pp.~334 and 337]{LLS}.

A related constant $\gamma_b'$, called the Berg\'e-Martinet constant,
is defined as follows:
\begin{equation}
\label{22'}
\gamma'_b = \sup\left\{ \lambda_1(L) \lambda_1(L^*)\left| L \subseteq
\R ^b \right. \right\},
\end{equation}
where the supremum is over all Euclidean lattices $L$.  Here $L^*$ is
the lattice dual to $L$.  If $L$ is the $\Z$-span of vectors $(x_i)$,
then $L^*$ is the $\Z$-span of a dual basis $(y_j)$ satisfying
$\langle x_i, y_j \rangle = \delta_{ij}$.

In a Riemannian manifold $(X,g)$, we define the volume
$\vol_k(\sigma)$ of a Lipschitz $k$-simplex $\sigma: \Delta^k
\rightarrow X$ to be the integral over the $k$-simplex~$\Delta^k$ of
the ``volume form'' of the pullback $\sigma^*(g)$.  The stable norm
$\|h\|$ of an element $h\in H_k(X,\R)$ is the infimum of the volumes
$\vol_k(c)=\Sigma_i |r_i| \vol_k(\sigma_i)$ over all real Lipschitz
cycles $c=\Sigma_i r_i \sigma_i$ representing $h$.  We define the
stable 1-systole of the metric $g$ by setting
\begin{equation}
\label{41b}
\stsys_1(g) =\lambda_1 \left( H_1^{\phantom{I}}(X, \Z)_\R, \|\cdot\|
\right) ,
\end{equation}
where $\|\cdot\|$ is the stable norm in homology associated with the
metric~$g$.

\section{A decomposition of the John ellipsoid}

The following statement may be known by convex set theorists.  A proof
may be found in in \cite{BI94}.  Recall that the John ellipsoid of a
convex set in Euclidean space is the unique ellipsoid of largest
volume inscribed in it \cite{MS}.

\begin{lemma} 
\label{ellipsoid-decomposition}
Let $(V^d,\|\cdot\|)$ be a Banach space.  Let $\|\cdot\|_E$ be the
Euclidean norm determined by the John ellipsoid of the unit ball of
$\|\cdot\|$. Then there exists a decomposition of $\|\cdot\|_E^2$ into
rank-1 quadratic forms:
$$
  \|\cdot\|_E^2 = \sum_{i=1}^N \lambda_i L_i^2
$$
such that $N\le \frac{d(d+1)}{2}+1$, $\lambda_i>0$ for all $i$,
$\sum\lambda_i=d$, and $L_i:V\to\R$ are linear functions with
$\|L_i\|^*=1$ where $\|\cdot\|^*$ is the dual norm to~$\|\cdot\|$.
\end{lemma}

\section{An area-nonexpanding map}

Let $X$ be a compact Riemannian manifold, $Y$ a topological space, and
let $\varphi:X\to Y$ be a continuous map inducing an epimorphism in
one-dimensional real homology.  Then one defines the {\it relative
stable norm} $\|\cdot\|_{st/\varphi}$ on $H_1(Y;\R)$ by
$$
  \|\alpha\|_{st/\varphi}
  = \inf\{\|\beta\|_{st}: \beta\in H_1(X;\R),\ \varphi_*(\beta)=\alpha \} .
$$
where $\|\cdot\|_{st}$ is the ordinary (``absolute'') stable norm.
The stable norm itself may be thought of as the relative stable norm
defined by the Abel-Jacobi map to the torus $H_1(X,\R)/ H_1 (X,
\Z)_{\R}$.

\begin{definition}
We will say that a Lipschitz map $\AJ: X \to M$ between Riemannian
manifolds is ``non-expanding on all $d$-dimensional areas'' if for
every smooth $d$-dimensional submanifold $Y$ of $X$, one has
$\vol_d(\AJ(Y))\le\vol_d(Y)$.
\end{definition}

Equivalently, $Jac(\AJ|_Y)\le 1$ wherever $\AJ|_Y$ is differentiable.

Let $X^n$ be a compact Riemannian manifold, $V^d$ a vector space and
$\Gamma$ a lattice in $V$.  We will identify $V$ and
$H_1(V/\Gamma;\R)$.

\begin{proposition} 
\label{area-nonexp}
Let $\varphi:X\to V/\Gamma$ be a continuous map inducing an
epimorphism of the fundamental groups and $\|\cdot\|_E$ denote the
Euclidean norm on $V$ defined by the John ellipsoid of the relative
stable norm~\hbox{$\|\cdot\|_{st/\varphi}$}.  Then there exists a
Lipschitz map $\AJ: X\to(V/\Gamma,\|\cdot\|_E)$ which is homotopic to
$\varphi$ and non-expanding on all $d$-dimensional areas, where
$d=\dim(V)$.
\end{proposition}


The proof of Propostion \ref{area-nonexp} appears at the end of this
section.  There is a natural isomorphism $\Gamma\simeq \pi_1(V/
\Gamma) \simeq \pi_1(X)/\ker(\varphi_*)$.  Consider a covering space
$\XX$ of $X$ defined by the subgroup $\ker(\varphi_*)$ of $\pi_1(X)$.
Then $\Gamma$ acts on $\XX$ as the deck group
$\pi_1(X)/\ker(\varphi_*)$.  This action will be written additively, as
in \eqref{41} below.  It is sufficient to construct a Lipschitz map
$\XX\to V$ which is $\Gamma$-equivariant and does not expand
$d$-dimensional areas.  We need the following lemma.

\begin{lemma} \label{lipfn}
For every linear function $L:V\to\R$ such that
$\|L\|_{st/\varphi}^*=1$ there exists a \short function
$f:\XX\to\R$ such that
\begin{equation}
\label{41}
f(x+v)=f(x)+L(v)
\end{equation}
for all $x\in\XX$ and $v\in\Gamma$.
\end{lemma}

\begin{proof}
Fix an $x_0\in\XX$, consider the orbit $\XX_0=\{x_0+v:v\in\Gamma\}$
and define a function $f_0:\XX_0\to\R$ by $f_0(x_0+v)=L(v)$.  Note
that $f_0$ satisfies \eqref{41} for $x\in\XX_0$.  For every
$v\in\Gamma$ and $x\in\XX$, one has $L(v)\le \|v\|_{st/\varphi}$ and
$\|v\|_{st/\varphi}$ is no greater than the distance between $x$ and
$x+v$. Hence $f_0$ is \short.  Every \short function defined on a
subset of a metric space admits a \short extension to the whole space,
by the triangle inequality.  Moreover, an extension can be chosen so
that the equivariance \eqref{41} is preserved.  For example, we can
set $f(x) = \inf \{ f_0(y)+|xy| : y\in\XX_0 \}$, where $|xy|$ denotes
the distance.
\end{proof}

\begin{proof}[Proof of Proposition \ref{area-nonexp}]
Applying Lemma \ref{ellipsoid-decomposition} to the norm $\|\cdot\|
_{st/ \varphi}$ yields a decomposition
\[
\|\cdot\|_E^2 = \sum_{i=1}^N \lambda_i L_i^2,
\]
where $\lambda_i>0$, $\sum \lambda_i=d$, $L_i\in V^*$ and
$\|L_i\|^*_{st/\varphi}=1$.  Then a linear map $L:V\to\R^N$ defined by
$$
  L(x) = (\lambda_1^{1/2} L_1(x), \lambda_2^{1/2} L_2(x), \dots,
  \lambda_N^{1/2} L_N(x)) 
$$
is an isometry from $(V,\|\cdot\|_E)$ onto a subspace $L(V)$ of
$\R^N$, equipped with the restriction of the standard coordinate
metric of $\R^N$.

By Lemma \ref{lipfn}, for every $i=1,2,\dots,N$ there exists a \short
function $f_i:\XX\to\R$ such that $f_i(x+v)=f_i(x)+L_i(v)$ for
all $x\in\XX$ and $v\in\Gamma$.  Define a map $F:\XX\to\R^N$
by
\begin{equation}
\label{42}
F(x) = (\lambda_1^{1/2} f_1(x), \lambda_2^{1/2} f_2(x), \dots,
\lambda_N^{1/2} f_N(x)) .
\end{equation}
Observe that both $L$ and $F$ are $\Gamma$-equivariant with respect to
the following action of $\Gamma$ on $\R^N$:
\[
\Gamma\times\R^N \to \R^N, \quad (v,x) \mapsto x+L(v).
\]
Now let $\Pr_{L(V)}:\R^N \to L(V)$ be the orthogonal projection to
$L(V)$.  Then the composition $L^{-1}\circ \Pr_{L(V)}\circ F$ is a
$\Gamma$-equivariant map from $\XX$ to~$V$.  Since the projection is
nonexpanding and $L$ is an isometry, it suffices to prove that the map
$F$ of \eqref{42} is nonexpanding on $d$-dimensional areas.

Let $Y$ be a smooth $d$-dimensional submanifold of $\XX$.  Since $F$
is Lipschitz, the restriction $F|_Y$ is differentiable a.e.\
on~$Y$. Let $y\in Y$ be such that $F|_Y$ is differentiable at $y$, and
let $A=d(F|_Y):T_yY\to \R^N$.  Then we obtain
$$
trace(A^*A) = \sum \lambda_i \left| d \left( f_i\mid_Y^{\phantom{I}}
\right) \right|^2 \le \sum \lambda_i = d
$$
since the functions $f_i$ are \short.  By the inequality of geometric
and arithmetic means, we have
$$
\begin{aligned}
Jac(F|_Y)(x) &= \det(A^*A)^{1/2} \\ & \le \left(\tfrac 1d
trace(A^*A)\right)^{d/2} \\ & \le 1 ,
\end{aligned}
$$
proving the proposition.
\end{proof}

\section{Proof of Theorem \ref{13} and Theorem \ref{21}}
\label{seven}

Consider the Jacobi torus $J_1(X)=H_1(X,\R)/H_1(X,\Z)_\R$ and the
Abel-Jacobi map $\AJ_X: X\to J_1(X)$ constructed in Proposition
\ref{area-nonexp}.

\begin{remark} 
\label{71a}
The map $\AJ_X$ can be replaced by a smooth one by an arbitrarily
small perturbation, in such a way as to expand $d$-dimensional areas
at most by a factor $1+\epsilon$.  Our main tool will be the coarea
formula, see below.  We can also avoid the above smoothing argument,
and use instead the current-theoretic version of the coarea formula,
relying on H. Federer's theory of ``slicing".  Given an integral
current $T$ in a smooth oriented manifold $M$ (\eg in our case $T =
[M]$), and a Lipschitz map $f:M \to N$, one can in a sense decompose
$T$ by $f$, obtaining currents (slices) $\langle T,f,y \rangle$
supported in the fibers $f^{-1}(y)$, for a.e. $y \in N$. If $T$ is a
cycle then so are the slices; if $T$ is a submanifold and $f$ is
smooth, then the slices at regular values of $f$ are just (integration
over) the typical fibers of the map.  The main properties of this
operation are given in \cite[Thm 4.3.2, p. 438]{Fe}.  In particular,
item 4.3.2(2) is a version of the coarea formula.  In section~4.4,
Federer shows that the usual homology groups of a reasonably good
space (\eg manifold in our case) coincide with the ones defined via
currents.
\end{remark}

\begin{proof}[Proof of Theorem \ref{13}] 
We exploit the coarea formula \cite[3.2.11]{Fe}, \cite[p.~267]{Ch2} as
in \cite[Theorem 7.5.B]{Gr1}.

Away from the negligible singular set, the smooth map is a submersion.
Therefore the metric on $X$ can be modified by a volume-preserving
deformation so that the map actually becomes {\em distance\/}
decreasing, up to an arbitrarily small amount.  Note, however, that
the coarea formula could be applied even without replacement by a
short map.  Formally we only need two facts:
\begin{enumerate}
\item
if the map is area-nonincreasing, then the volume is no smaller than
the area of the image times the minimal area of a fiber,
\item
almost every fiber is ``typical'' and hence has area no less than the
generalized degree $\deg(\AJ_X)$, \cf inequality \eqref{71b}.
\end{enumerate}

Let $S=\AJ_X^{-1}(p)$ be the surface which is a regular fiber of least
$(n-b)$-volume.  (If there is none, choose $S$ to be within
$\epsilon>0$ of the infimum, and then let $\epsilon\to 0$.)  Then
\begin{equation}
\label{71b}
\deg(\AJ_X) \stsys_1(\g)^b \le \vol_{n-b}(S) \, \stsys_1(\g)^b .
\end{equation}
Note that $\AJ_X$ induces isometry in $H_1(\;,\R)$ with respect to the
stable norm of the metric $g$.  Hence
\[
\deg(\AJ_X) \stsys_1(\g)^b \leq \vol_{n-b}(S) \, \stsys_1 (J_1(X),
\|\cdot\|)^b.
\]
We now replace the stable norm $\|\cdot\|$ by the flat Euclidean
metric $\|\cdot\|_E$ defined by the John ellipsoid of the stable norm:
\[
\deg(\AJ_X) \stsys_1(\g)^b \leq \vol_{n-b}(S) \;
\stsys_1(J_1(X),\|\cdot\|_E)^b .
\]
By definition of the Hermite constant,
\[
\deg(\AJ_X) \stsys_1(\g)^b \leq \vol_{n-b}(S)\; \gamma_b^{\frac{b}{2}}
\vol_b(J_1(X), \|\cdot\|_E).
\]
Now we apply the coarea formula to our map which is decreasing on
$b$-dimensional volumes, to obtain $\vol_{n-b}(S) \vol_b(J_1(X),
\|\cdot\|_E)\leq \vol_n(X)$, completing the proof.
\end{proof}

\begin{proof}[Proof of Theorem~\ref{21}] 
We apply Theorem~\ref{13} together with the inequality $\pisys_1^2(g)<
\sigma_2 \deg(\AJ_X)$.  Here Pu's inequality does not suffice.
Indeed, a typical fiber of $\AJ_X$ may not be diffeomorphic to $\R
P^2$.  An application of Pu's inequality \eqref{pu3} yields a suitably
short loop which is essential in the typical fiber.  However, a loop
which is essential in the typical fiber, may not be essential in the
ambient manifold $X$.  Thus, we need a generalisation of Pu's
inequality.  The required generalisation is inequality~\eqref{23}
above, \cf \cite{CK}, applied to the composed map $\hat \phi: S\to \XX
\to K$.
\end{proof}

\section*{Acknowledgments}
We are grateful to J. Fu for help with integral currents and slices in
section \ref{seven}.

 \end{document}